\documentclass{article}

\usepackage{hyperref} 
\usepackage[super, sort]{natbib}
\bibpunct{(}{)}{;}{a}{,}{,} 
\usepackage{amsmath}
\usepackage{url}

\begin{document}

\title{Fermi-Dirac Integrals in terms of Zeta Functions }
\author{Michael Morales\footnote{Currently member of Department of Mathematics, Ohio University, OH 45701, USA}\\
 \small{Department of Physics, Universidad del Valle de Guatemala,}\\ \small{ Guatemala City 01015, Guatemala. }}
\date{}
\maketitle

\begin{abstract}

A constructive method to write Fermi-Dirac integrals in terms of Riemann and Hurwitz zeta functions is presented. By using of an auxiliary function we can decompose Fermi-Dirac integrals in terms of well-known special functions, which evades the necessity of iterative methods. The procedure is presented for a particular example; neverthless, its extension to general Fermi-Dirac integrals is straightforward. The final expression obtained along with Sommerfeld's lemma gives nearly a complete representation of Fermi-Dirac integrals.

\end{abstract}

\section{Introduction}

 Fermi-Dirac integrals $F_\frac{k}{2}(\eta)$  are defined as

\begin{center} 
\begin{equation}F_{\frac{k}{2}}(\eta)=\displaystyle\int_{0}^{\infty} \frac{\xi^\frac{k}{2}}{1+e^{\xi-\eta}}d\xi\end{equation}
\label{eqn:fdgral}
\end{center}

appear in a variety of physical disciplines including partial degenerate stars \citep{clayton}, thermal conductions by electrons \citep{khal}, and condensed matter physics in general, where eta = mu/kT, k is the Boltzmann constant, T is the temperature and mu the chemical potential. In this work we determine expressions that allow us to calculate Fermi-Dirac integrals specifically for eta < 5, values that so far have been calculated using iterative methods as the trapezoidal method \citep{mohankumar}. In the following sections we present the case k=1; nonetheless, the presented procedure can be extended to any integer k.\\

This procedure involves a basic auxiliary function that relates $F_\frac{1}{2}(\eta)$ with $F_\frac{1}{2}(0)$, which in practical terms is a model based on their graphical behavior. This makes the integral $F_\frac{1}{2}(\eta)$ easier to calculate obtaining Riemann and Hurwitz Zeta functions, with almost no need for iterative methods.

\section{Sommerfeld's lemma applied to $F_\frac{1}{2}(\eta)$}

       Let $\psi (u)$ be a sufficiently regular function which vanishes for $u=0$, then we have the asymptotic formula \citep{chandra}  \\

$ \displaystyle\int_0^\infty \frac{du}{\frac{1}{\Lambda}e^u+1}\frac{d\psi (u)}{du} =$
\begin{equation} \psi (u_0) +2\left[ c_2 \psi ''(u_o)+c_4 \psi ^{(iv)} (u_o)+\cdots \right]
\end{equation}
  
\noindent where $u_o=\log[ \Lambda ]$ and $c_2, c_4, \ldots$ are numerical coefficients defined by

\begin{equation}
 c_\nu = 1-\frac{1}{2^\nu}+\frac{1}{3^\nu}-\frac{1}{4^\nu}+\cdots
\end{equation}

\noindent It can be show that the Lemma applied to $F_\frac{1}{2}(\eta)$ in terms of Riemann Zeta function $\zeta(z)$, with $\psi(u)=u^{1/2}$, and $u_0=\eta$, results in

$
\displaystyle\int_0^\infty \frac{u^{1/2}}{1+e^{u-\eta}}du=$
\begin{equation}
\frac{2}{3}\eta^{3/2}+2\displaystyle\sum_{n=1}^{\infty} (1-2^{1-2n})\zeta (2n)\frac{(4n-5)!!}{2^{2n-1}} \eta^{-\frac{4n-3}{2}}
\label{eqn:sommerzeta}
\end{equation}

\noindent which is valid when we neglect values of order $e^{-\eta}$. In fact, it can be shown that (\ref{eqn:sommerzeta})  is valid for $\eta >3$.

\section{$F_\frac{1}{2}(\eta)$ and the auxiliary function}
 
  As mentioned, we analyze the Fermi-Dirac integral $F_\frac{1}{2}(\eta)$ 

\begin{center} 
\begin{equation}F_\frac{1}{2}(\eta)=\displaystyle\int_{0}^{\infty} \frac{\xi^\frac{1}{2}}{1+e^{\xi-\eta}}d\xi
\label{eqn:FM12}
\end{equation}
\end{center}

\noindent by  focusing on the integrand

\begin{center} 
\begin{equation}f_\frac{1}{2}(\xi,\eta)= \frac{\xi^\frac{1}{2}}{1+e^{\xi-\eta}}
\label{eqn:fm12}
\end{equation} 
\end{center}

\noindent and the related function $f_\frac{1}{2}(\xi)$

\begin{center} 
\begin{equation}f_\frac{1}{2}(\xi)= \frac{\xi^\frac{1}{2}}{1+e^{\xi}}
\label{eqn:fm0}
\end{equation} 
\end{center}

\noindent  We also define a function, denoted as $f(\xi,\eta)$, as the ratio of both functions, $f_\frac{1}{2}(\xi,\eta)$ and  $f_\frac{1}{2}(0)$

\begin{center} 
\begin{equation}f(\xi,\eta)=\frac{f_\frac{1}{2}(\eta)}{f_\frac{1}{2}(0)}
                                      =\frac{1+e^{\xi}}{1+e^{\xi-\eta}}
\label{eqn:ratio}
\end{equation} 
\end{center}  

\noindent this function, graphically, behaves as $a-b e^{-c \xi}$. This suggests that we can write (\ref:{eqn:ratio}) as

\begin{center} 
\begin{equation}f(\xi,\eta)=a-b e^{-c \xi}
\label{eqn:f1}
\end{equation} 
\end{center}  

\noindent where $a, b,$ and $c$ are constants. The first two can be found considering some analytical aspects of equation (\ref{eqn:ratio}):

\begin{center} 
\begin{equation}
f(0,\eta)=\frac{2}{1+e^{-\eta}}
\end{equation}
\end{center}  

\noindent and 

\begin{center} 
\begin{equation}
\displaystyle\lim_{x\to\infty}f(0,\eta)=e^{\eta}
\end{equation}
\end{center}    

\noindent this implies that $a=e^\eta$ and $b=\frac{e^\eta-1}{e^{-\eta}+1}$. With this we obtain 

\begin{center} 
\begin{equation}f(\xi,\eta)=e^\eta-\frac{e^\eta-1}{e^{-\eta}+1} e^{-c \xi}
\label{eqn:f2}
\end{equation} 
\end{center}  

\noindent This leaves the issue of finding  $c$. We take the following approach, perhaps the only one in this procedure that needs iterative methods: we find the extreme point of $f_\frac{1}{2}(\xi,\eta)$, equation (\ref{eqn:fm12}). After some algebra, we obtain the following equation which must solved for $\xi$

\begin{center}
\begin{equation}
e^\eta+e^\xi(1-2\xi)=0
\label{eqn:xmax}
\end{equation}
\end{center}

\noindent We call the solution of (\ref{eqn:xmax}) $\xi_m$, which depends on $\eta$. Using equations (\ref{eqn:ratio}) and (\ref{eqn:f2}), and substituting $\xi_m$, we obtain:

\begin{center}
\begin{equation}
c=-\frac{1}{\xi_m}Ln\left[  \frac{(e^{-\eta}+1)\left(e^\eta-\frac{1+e^{\xi_m}}{1+e^{\xi_m-\eta}}\right)}{e^\eta-1}     \right]
\label{eqn:c}
\end{equation}
\end{center}

\noindent  It can be shown that calculating the maxima of equation (\ref{eqn:fm12}) and the product $f(\xi,\eta)  f_\frac{1}{2}(\xi)$ leads to the same result. After we determine the constants in equation (\ref{eqn:f1}), we can write equation (\ref{eqn:FM12}) as 
 
\begin{center} 
\begin{eqnarray}
\displaystyle\int_{0}^{\infty} \frac{\xi^\frac{1}{2}}{1+e^{\xi-\eta}}d\xi=\displaystyle\int_{0}^{\infty} \left(e^\eta-\frac{e^\eta-1}{e^{-\eta}+1} e^{-c \xi}\right)\frac{\xi^\frac{1}{2}}{1+e^{\xi}}d\xi
=e^\eta\displaystyle\int_{0}^{\infty}\frac{\xi^\frac{1}{2}}{1+e^{\xi}}d\xi-\frac{e^\eta-1}{e^{-\eta}+1}\displaystyle\int_{0}^{\infty}\frac{e^{-c \xi}\xi^\frac{1}{2}}{1+e^{\xi}}d\xi
\label{eqn:preintegral}
\end{eqnarray}
\end{center}

\noindent where the first integral in the right hand side can be expressed as \citet{arfken} 

\begin{center}
\begin{equation}
\displaystyle\int_{0}^{\infty}\frac{\xi^\frac{1}{2}}{1+e^{\xi}}d\xi=\Gamma\left(\frac{3}{2}\right) \left(1-2^{-\frac{1}{2}}\right)\zeta\left(1+\frac{1}{2}\right)
\end{equation}
\end{center}

 The second term in (ref{eqn:preintegral}) after expanding some terms takes the form

\begin{center}
\begin{eqnarray}
\displaystyle\int_{0}^{\infty}\frac{e^{-c \xi}\xi^\frac{1}{2}}{1+e^{\xi}}d\xi=\Gamma\left(\frac{3}{2}\right)\displaystyle\sum_{n=0}^{\infty}\frac{(-1)^n}{(c+1+n)^{1+\frac{1}{2}}}
=\Gamma\left(\frac{3}{2}\right)\left[ 2^{-\frac{1}{2}}\zeta\left(1+\frac{1}{2},\frac{c+1}{2}\right)-\zeta\left(1+\frac{1}{2},c+1\right)\right]
\end{eqnarray}
\end{center}

In this expression we recognize $\zeta(p,q)$ as the Hurwitz zeta function, where we have used the result in \citet{williams}. Using these results we can finally write $F_\frac{1}{2}(\eta)$ as

\begin{center}
\begin{equation}
F_\frac{1}{2}(\eta)=\Gamma\left(\frac{3}{2}\right)\left[e^\eta\left(1-2^{-\frac{1}{2}}\right)\zeta\left(1+\frac{1}{2}\right)-\frac{e^\eta-1}{e^{-\eta}+1}\left[ 2^{-\frac{1}{2}}\zeta\left(1+\frac{1}{2},\frac{c+1}{2}\right)-\zeta\left(1+\frac{1}{2},c+1\right)\right]\right]
\label{eqn:final}
\end{equation}
\end{center}

\section{Some results}

  Using equation (\ref{eqn:final}) we construct Table 1 for some specific values of $\eta$ which allows a comparison with values taken from \citet{mcdougall}, including some values for positive $\eta<1$. These values have been used in electron conduction opacity, \citep{khal}, showing how reliable and potentially applicable this method is. \\

\begin{table}[hbt]
\caption{ $ F_\frac{1}{2}(\eta) $ values comparison}
\centering
\begin{tabular}{c c c c }
\hline
$\eta$&$F_\frac{1}{2}$&$F_\frac{1}{2}$&difference (\%)\\
           & eq. (\ref{eqn:final})&  McDougall \& Stoner &\\

 \hline
\hline
-4.0 &  0.01614  & 0.01613 &  0.06\\
-3.0 &  0.04345  & 0.04337 &  0.18\\
-2.0 &  0.11506  & 0.11459  & 0.41\\
-1.0 &  0.29241  & 0.29050  & 0.66\\
0.0 &   0.67809  & 0.67810  & 0.00\\
0.1 &   0.73203  & 0.73340  & 0.19\\
0.5 &   0.97795  & 0.99021  & 1.24\\
1.0 &   1.35129  & 1.39638  & 3.23\\
2.0 &   2.30003     & 2.50246  & 8.09\\
3.0 &   3.58315     & 3.97699  & 9.90\\
4.0 &   5.54950       & 5.77073  & 3.83\\
5.0 &   8.99919     & 7.83797    & 14.82\\

\hline
\end{tabular}
\end{table}

\section {Observations}
 

  The model according to (\ref{eqn:f1}) is reliable for $\eta < 5$, specially with for values close to zero and negative $\eta$. Yet, as the ratio of the functions in equation (\ref{eqn:ratio}) shows a higher inflexion point for $\eta>5$, the model needs an improvement. The reason for the limitation around this value appears when analyzing the equation (\ref{eqn:final}), where a simple examination shows a dependency of order $e^\eta$; recall that, in Sommerfeld lemma, this dependence is of order $e^{-\eta}$.\\

  The use of the auxiliary function, a simple exponential model, practically evades the use of iterative methods or pole expansion methods. Numerical methods are needed only for the simple equation (\ref{eqn:xmax}) and to determine Riemann and Hurwitz zeta functions. 

\section{Conclusion}

  Equation (\ref{eqn:final}) gives an expression for calculating  Fermi-Dirac integrals $F_\frac{1}{2}(\eta)$ in terms of the Riemann and Hurwitz Zeta functions. This approach complements Sommerfeld's lemma, resulting in a nearly complete expression for Fermi-Dirac integrals. \\

\noindent Using the same procedure we can find expressions the corresponding expression $F_\frac{3}{2}(\eta)$

\begin{center}
\begin{eqnarray}
&& F_\frac{3}{2}(\eta)=\Gamma\left(\frac{5}{2}\right)\left[e^\eta\left(1-2^{-\frac{3}{2}}\right)\zeta\left(1+\frac{3}{2}\right) \right. \nonumber \\ 
&& \left.-\frac{e^\eta-1}{e^{-\eta}+1}\left[ 2^{-\frac{3}{2}}\zeta\left(1+\frac{3}{2},\frac{c+1}{2}\right)-\zeta\left(1+\frac{3}{2},c+1\right)\right]\right]
\end{eqnarray}
\end{center}

 \noindent in general, the procedure presented here can be extended for any integer $k$ to determine $F_\frac{k}{2}(\eta)$,

\begin{center}
\begin{eqnarray}
&& F_\frac{k}{2}(\eta)=\Gamma\left(1+\frac{k}{2}\right)\left[e^\eta\left(1-2^{-\frac{k}{2}}\right)\zeta\left(1+\frac{k}{2}\right)\right. \nonumber \\
&& \left.-\frac{e^\eta-1}{e^{-\eta}+1}\left[ 2^{-\frac{k}{2}}\zeta\left(1+\frac{k}{2},\frac{c+1}{2}\right)-\zeta\left(1+\frac{k}{2},c+1\right)\right]\right]
\end{eqnarray}
\end{center}

\noindent with the remark that $\xi_m$, in order to find $c$ (notice that $c$ and $f(\xi,\eta)$ are $k$-independent) equation (\ref{eqn:xmax})  must also be generalized:

\begin{center}
\begin{equation}
k e^\eta+e^\xi(k-2\xi)=0
\end{equation}
\end{center}

  Although is not a complete close form expression, by using an elementary methods, we obtained a single expression with which when combined with Sommerfeld's lemma, we have a complete accurate description for Fermi-Dirac integrals $F_{\frac{k}{2}}(\eta)$.

\bibliographystyle{plainnat} 
\bibpunct{(}{)}{;}{a}{,}{,} 

\bigskip



\bigskip

\end{document}